\documentclass[11pt, oneside]{article}   	
\usepackage{geometry}                		
\usepackage{graphicx}				
\newtheorem{thm}{Theorem}
\newtheorem{lemma}[thm]{Lemma}
\newtheorem{conj}[thm]{Conjecture}

\usepackage{amssymb}

\title{Hilbert series of generic ideals in products of projective spaces}
\author{Ralf Fr\"oberg}
\date{}							

\begin{document}
\maketitle

\begin{abstract}
If $k[x_1,\ldots,x_n]/I=R=\sum_{i\ge0}R_i$, $k$ a field, is a standard graded algebra, the 
Hilbert series of $R$ 
is the formal power series $\sum_{i\ge0}\dim_kR_it^i$. It is known already since Macaulay
which power series are Hilbert 
series of graded algebras \cite{ma}. A much harder question is which series are Hilbert series if
we fix the number of generators of $I$ and their degrees, say for ideals $I=(f_1,\ldots,f_r)$, 
$\deg f_i=d_i$, $i=1,\ldots,r$. In some sense "most" ideals with fixed degrees of their generators
have the same Hilbert series. There is a conjecture for the Hilbert series of those "generic"
ideals, see below. In this paper we make a conjecture, and
prove it in some cases, in the case of generic 
ideals of fixed degrees in the coordinate ring of ${\mathbb P}^1\times{\mathbb P}^1$, which might
be easier to prove.


\end{abstract}

\section{Background}
The conjecture and the results in section 2 below are inspired by the corresponding conjecture 
and results in the singly graded case.

\begin{conj}[\cite{fr}]
Let $I=(f_1,\ldots,f_r)\subset k[x_1,\ldots,x_n]$, $k$ an infinite field,
be an ideal generated by generic forms $f_i$
with $\deg(f_i)=d_i$, and let
$R= k[x_1,\ldots,x_n]/I$. Then
$$R(u)=[\prod_{i=1}^r(1-u^{d_i})/(1-u)^n]_+$$
here $[\sum a_iu^i]_+=\sum b_iu^i$, where $b_i=a_i$ if for all $j\le i$ we have $a_j>0$, 
and $b_i=0$ otherwise.
\end{conj}
We first comment on the use of the word "generic". A polynomial of degree $d$ in
$k[x_1,\ldots,x_n]$ is a linear combination of ${n+d-1\choose d}$ monomials. Thus an ideal
$(f_1,\ldots,f_r)$, $\deg(f_i)=d_i$, can be considered as a point in $A={\mathbb A}^N$,
$N=\sum_{i=1}^r{n+d_i-1\choose d_i}$. There is a Zariski-open subset of $A$, for which
the Hilbert series is constant. Ideals corresponding to points in that Zariski-open set are
what we call generic, see \cite{fr-lo}.

\bigskip
The conjecture is proved for $r\le n$ (trivial), for $n\le2$ \cite{fr}, for $n=3$ \cite{an}, 
for $r=n+1$ \cite{st}. There
are partial results in \cite{ho-la}, \cite{fr-ho}, \cite{au}, \cite{mi-mi}, \cite{ba-on}, \cite{ni}, \cite{ne}.

\bigskip
If $J=(l_1^{d_1},\ldots,l_r^{d_r})$, where $l_i$ are generic linear forms. Sometimes, but not
always, the Hilbert series of $k[x_1,\ldots,x_n]/J$ equals the one in the conjecture. There
is a conjecture on when it does, \cite{ia,ch}.

\section{${\mathbb P}^1\times{\mathbb P}^1$}
 We are considering homogeneous ideals in the coordinate ring of 
 ${\mathbb P}^1\times{\mathbb P}^1$. Thus, let $k$ be an infinite field,
$S=k[x_0,,x_1,y_0,y_1]$ be bigraded, $\deg(x_i)=(1,0)$, $\deg(y_j)=(0,1)$, and let  $I$ 
be a bihomogeneous ideal, so generated  by bihomogeneous elements.
Hence $R=S/I$ is bigraded, $R=\oplus_{i,j\ge0}R_{i,j}$. The Hilbert series  of $R$ is defined as
$R(u,v)=\sum\dim_kR_{i,j}u^iv^j$. 
We are  interested in the case when the ideal is generated by "generic" elements.
Given a sequence
of degrees ${\bf d}_1,\ldots,{\bf d}_r$, we denote the space of ideals $I=(f_1,\ldots,f_r)$ where 
$\deg f_i={\bf d}_i$ by ${\bf I}_{{\bf d}_1,\ldots,{\bf d}_r}$. 
An element of degree $(d,e)$ is a linear combination of
$(d+1)(e+1)$ monomials. Thus an ideal in  ${\bf I}_{{\bf d}_1,\ldots,{\bf d}_r}$ can be
considered as a point in ${\mathbb A}_k^N$, $N=\sum_{i=1}^r(d_i+1)(e_i+1)$ where 
${\bf d}_i=(d_i,e_i)$. We partially order Hilbert series termwise, so that 
$\sum a_{ij}u^iv^j\ge\sum b_{ij}u^iv^j$ if $a_{ij}\ge b_{ij}$ for all $i,j$.

\begin{thm}
There are only a finite number of possibilities for Hilbert series of ideals in  
${\bf I}_{{\bf d}_1,\ldots,{\bf d}_r}$.
There is a nonempty Zariski open 
part of  of ${\bf I}_{{\bf d}_1,\ldots,{\bf d}_r}$ where the Hilbert series
is constant. This constant Hilbert series is the smallest possible for ideals in
 ${\bf I}_{{\bf d}_1,\ldots,{\bf d}_r}$.
 \end {thm}
 {\bf Proof} The corresponding theorems in the singly graded case, \cite[Theorem 1]{fr-lo} ,
 \cite[Theorem p.120]{fr}, and \cite[Proposition 1]{fr-gu-lo}  are easily adapted. We call points in this 
 nonempty Zariski open set generic.
 
 \medskip
 We define $(d,e)\le(f,g)$ if $d\le f$ and $e\le g$, and $(d,e)\ge(f,g)$ if $d\ge f$ and $e\ge g$.
 Furthermore $[\sum_{i,j} a_{ij}u^iv^j]_+=\sum_{i,j}b_{ij}u^iv^j$, where $b_{ij}=a_{ij}$ if
 $a_{kl}>0$ for all $(k,l)\le(i,j)$ and $b_{i,j}=0$ otherwise.
 
 \begin{lemma}
 Let $R=S/I$ be bigraded and $f\in R_{i,j}$. Then $(R/f)(u,v)\ge[(1-t^iu^j)(R(u,v)]_+$,
 \end{lemma}
{\bf Proof} Consider the map $f\cdot\colon R_{d-i,e-j}\rightarrow R_{d,e}$. The image is largest
 if the map is of maximal rank, i.e., either injective or surjective, so 
 $\dim_k(R/f)_{d,e}\ge\max\{0,\dim R_{d,e}-\dim R_{i,j}\}$. If $\dim (R/f)_{d,e}=0$, then
 $\dim (R/f)_{d+f,e+g}=0$ for all $(f,g)\ge0$.
 
 \begin{lemma}
 $[(1-u^iv^j)[\sum_{i,j}a_{ij}u^iv^j]_+]_+=[(1-u^iv^j)\sum_{i,j}a_{ij}u^iv^j]_+$.
 \end{lemma}
 {\bf Proof} Easy calculation.
 
 \medskip
 These two lemmas give the following.
 \begin{thm}
 Let $I=(f_1,\ldots,f_r)\subset k[x_0,x_1,y_0,y_1]=S$, $\deg f_i=(d_i,e_i)$. Then
 $S/I(u,v)\ge[\prod_{i=1}^r(1-u^{d_i}v^{e_i})/((1-u)^2(1-v)^2)]_+$.
 \end{thm}
 
 We now give a conjecture in the case when the $f_i$'s are generic,  c.f. \cite{fr-lu}. To prove the conjecture
 for some fixed $({\bf d}_1,\ldots,{\bf d}_r)$ it suffices to give one example with the conjectured
 series. If the conjecture is true for these parameters, then almost all ideals have the conjectured
 series, so we must be very unlucky if we miss the series with a random choice of coefficients.
 
 \begin{conj}
 Let $I=(f_1,\ldots,f_r)\subset k[x_0,x_1,y_0,y_1]=S$, $\deg f_i=(d_i,e_i)$ generic. Then
 $(S/I)(u,v)=[\prod_{i=1}^r(1-u^{d_i}v^{e_i})/((1-u)^2(1-v)^2)]_+$.
 \end{conj}
 
 We have checked that the conjecture is true in the following cases. Some of these were checked 
 by Alessandro Oneto. Except for the first class, we have used computer calculations.
 
 \smallskip\noindent
 1. For small $r$ the concepts of ideal generated by generic forms and complete intersection
 agrees. It is well known that the conjecture is true for complete intersections.

 \smallskip\noindent
 2. $\deg f_i=(1,1)$ for all $i$, any $r$.
 
 \smallskip\noindent
 3. Some $f_i$ of degree (1,1), some of degree (1,2), any $r$.
 
 \smallskip\noindent
 4. $\deg f_i=(1,2)$ for all $i$, any $r$.
 
 \smallskip\noindent
5. $\deg f_i=(2,2)$ for all $i$, any $r$.
 
 \medskip
 We also checked that the corresponing conjecture is true for $\deg f_i=(1,1,1)$ for all $i$, any $r$, in
 ${\mathbb P}^1\times{\mathbb P}^1\times{\mathbb P}^1$.
 
 \medskip
 On the other hand, the corresponding conjecture in ${\mathbb P}^2\times{\mathbb P}^2$ cannot
 be true. For four generic forms of degree (1,1) the conjecture would give that $R_{d,d}=0$
  if $d>>0$. The correct statement is that $\dim_kR_{d,d}=6$ if $d>>0$.
  
  \medskip
  We think that the conjecture is challenging enough, but we also give some questions.
  
  \medskip\noindent
  {\bf Question} What is the Hilbertseries for generic ideals in 
  $( {\mathbb P}^1)^k={\mathbb P}^1\times\cdots\times{\mathbb P}^1$, $k>2$ times?
  
  \medskip\noindent
  {\bf Question}  What is the Hilbertseries for generic ideals in 
  ${\mathbb P}^m\times{\mathbb P}^n$?
  
  \medskip\noindent
  {\bf Question} Let $f_i$ be generic linear forms in $k[x_1,x_2]$ and $g_i$
  generic linear forms in $k[y_1,y_2]$, and let $I=(f_1^{d_1}g_1^{e_1},\ldots,f_r^{d_r}g_r^{e_r})$.
  What is the Hilbert series of $k[x_1,x_2,y_1,y_2]/I$?

\end{document}